\documentclass[12pt]{article}
\usepackage[utf8]{inputenc}
\usepackage{amssymb,amsmath,amsfonts,eucal,mathrsfs,amsthm} 
\usepackage[colorinlistoftodos]{todonotes}
\usepackage{enumitem}
\newtheorem{theorem}{Theorem}
\newtheorem{proposition}[theorem]{Proposition}

\newtheorem{remarks}[theorem]{Remarks}

\theoremstyle{definition}
\newcommand{\R}{\mathbb{R}}

\newcommand{\Sf}{\mathbb{S}}

\newcommand{\Hy}{\mathbb{H}}

\newcommand{\spa}{\mbox{span}}
\newcommand{\hess}{\mbox{Hess\,}}
\newcommand{\Ric}{\mbox{Ric}}

\def\<{{\langle}}
\def\>{{\rangle}}

\def\n{\nabla}

\def\be{\begin{equation} }
\def\ee{\end{equation} }

\begin{document}

\title{On Einstein submanifolds 
of Euclidean space}
\author{M. Dajczer, C.-R. Onti and Th. Vlachos}
\date{}
\maketitle

\begin{abstract}
Let the warped product  $M^n=L^m\times_\varphi F^{n-m}$, 
$n\geq m+3\geq 8$, of Riemannian manifolds be an Einstein 
manifold with Ricci  curvature $\rho$ that admits an isometric 
immersion into Euclidean space with codimension two. Under the 
assumption that $L^m$ is also Einstein, but not of constant 
sectional curvature, it is shown that  $\rho=0$ and that the 
submanifold is locally a cylinder with an Euclidean factor of 
dimension at least $n-m$. Hence $L^m$ is also Ricci flat.
If $M^n$ is complete, then the same conclusion holds globally 
if the assumption on $L^m$ is replaced by the much weaker 
condition that either its scalar curvature $S_L$ is constant 
or that $S_L\leq (2m-n)\rho$.
\end{abstract}

A Riemannian manifold $M^n$ is said to be \emph{Einstein} if 
its Ricci tensor is proportional to the metric, that is, if 
$\Ric_M(X,Y)=\rho\,\<X,Y\>$ for any vector fields 
$X,Y\in\mathfrak{X}(M)$ and a constant $\rho\in\R$. 
Then $\rho$ is the not normalized constant Ricci curvature 
of $M^n$.  Einstein manifolds produced as warped product have 
been quite intensively studied since Besse \cite{Be}, in his 
book, presented many results on this subject. On the 
contrary, this class of  manifolds has seldom been considered 
as submanifolds of Euclidean space.

Cartan \cite{Th} and Fialkow \cite{Fi} proved that any Einstein 
hypersurface in Euclidean space $f\colon M^n\to\R^{n+1}$, $n\geq 3$, 
has constant sectional curvature, namely, it is either 
flat or an open subset of a round sphere. In \cite{DOV2} we 
investigated the case of Einstein submanifolds $M^n$ in Euclidean 
space $\R^{n+2}$ when $M^n$ possess the metric structure of a 
warped product over a surface, that is, if
$M^n=L^2\times_\varphi F^{n-2}$. We constructed many local 
examples with non-constant sectional curvature. On the other
hand, we showed that if $F^{n-2}$ has constant sectional 
curvature, the only complete examples with dimension 
$n\geq 5$ are a certain product of round spheres and the 
Ricci flat Generalized Schwarzschild metric (GS-metric) 
isometrically immersed as an $(n-2)$-rotational submanifold. 

Recall that the GS-metric is a Ricci flat warped product 
metric on the product manifold $\R^2\times F^{n-2}$ where
the constant Ricci curvature of $F^{n-2}$ is positive. The 
warped product $M^n=\R^2\times_\varphi\Sf^{n-2}(1)$ is said 
to be endowed with the GS-metric metric when $\R^2$ in polar 
coordinates has the rotational invariant warped metric 
$ds^2=dt^2+\varphi'^2(t)d\theta^2$ where the warping function
$\varphi\in C^\infty([0,+\infty))$ is the unique positive 
solution of the differential equation $\varphi'^2=1+c/\varphi^{n-3}$ 
with a given initial condition  and a constant $c<0$. 
\vspace{1ex}

The above considerations raise the following general question: 
\emph{Which are the Einstein submanifolds 
$f\colon M^n=L^m\times_\varphi F^{n-m}\to\R^{n+2}$ when 
$m\geq 3$ and $n-m\geq 2$?}  
\vspace{1ex}

In this paper we investigate the case of isometric 
immersions of Einstein manifolds of the form 
$M^n=L^m\times_\varphi F^{n-m}$ when the base $L^m$ has
dimension $m\geq 5$. We recall that if $M^n$ is Einstein, 
then  $F^{n-m}$ also has to be Einstein; cf.\ Corollary $9.107$ 
in \cite{Be}. In addition, we have that $M^n$ is a complete 
manifold if and only if both factors are complete. 

As already clarified in \cite{DOV2}, in order to obtain classification 
results, the assumption that an Euclidean submanifold is Einstein 
is quite weak, thus requiring the use of additional hypothesis in 
order to be successful.
In the local case, we assume that $L^m$ is also an Einstein manifold. 
This class of manifolds have already been intrinsically under 
consideration, for instance in \cite{BG}. In the case that the
manifold is complete, we do with a quite weaker assumptions on the
scalar curvature of $L^m$, for instance, that it is constant.
Intrinsically, the Einstein manifolds satisfying that 
condition have been considered in \cite{HPW}.\vspace{1ex}

We call \emph{nontrivial} an Einstein manifold $N^n$, $n\geq 4$,  
if in no open subset the sectional curvature is a nonnegative 
constant. Due to the aforementioned result by Fialkow, the 
condition of non triviality implies that non open subset of 
$N^n$ admits an isometric immersion into $\R^{n+1}$.
\vspace{1ex}

Locally, and when the base is Einstein, we have the following result.

\begin{theorem}\label{local}
Let $M^n=L^m\times_\varphi F^{n-m}$, $n\geq m+3\geq 8$, and $L^m$
be nontrivial Einstein manifolds. If there exists an isometric 
immersion $f\colon M^n\to\R^{n+2}$ then $M^n$ is Ricci flat and 
the warping function satisfies that $\|\n\varphi\|=c\geq 0$ is 
constant. Moreover, any point of an open dense subset $M^n_0$ 
of $M^n$ lies in an open product neighborhood 
$U=L^m_0\times F^{n-m}_0$, where 
$L^m_0\subset L^m$ and $F^{n-m}_0\subset F^{n-m}$, such 
that on $M^n_0$  always the same of one of the following holds:
\begin{itemize}
\item[(i)] We have  $c=0$, there is an isometric immersion 
$f_0\colon L^m_0\to\R^{m+2}$, 
$F^{n-m}_0\subset\R^{n-m}$ is an open subset and 
$f|_U=f_0\times id_{F_0}$.
\item[(ii)]  
We have $c>0$, $L_0^m=I\times N^{m-1}$ where $I\subset\R$ is 
an open interval and  $F_0^{n-m}\subset\Sf^{n-m}(1/\sqrt{c})$ 
is an open subset. There is an orthogonal splitting 
$\R^{n+2}=\R^{m+1}\oplus\R^{n-m+1}$ such that $f|_U=f_0\times j$ 
where $f_0\colon N^{m-1}\to\R^{m+1}$ is an isometric immersion 
and the map $j\colon I\times F_0^{n-m}\to\R^{n-m+1}$ given by 
$j(t,y)=ty$ is a parametrization of an open subset of $\R^{n-m+1}$.
\end{itemize}
\end{theorem}

\vspace{3ex}

\begin{picture}(150,84)\hspace{-15ex}
\put(90,30){$I\times\Sf^{n-m}(1/\sqrt{c})
\supset(I\times F_0^{n-m})\times N^{m-1}=U\subset M^n$}
\put(235,55){$j$} 
\put(274,55){$f_0$}
\put(230,42){\vector(0,1){30}} \put(270,42){\vector(0,1){30}}
\put(206,80){$\R^{n-m+1}\oplus\R^{m+1}=\R^{n+2}$}
\end{picture}
\vspace{-2ex}

In case $(ii)$ we have that $U$ is foliated by hypersurfaces
parametrized by $I$ whose images are $(n-m)$-rotational 
submanifolds. Besides, observe that the submanifold $f(U)$ in 
part $(ii)$ is an $(n-m+1)$-cylinder and thus it is trivially 
an $(n-m)$-cylinder as the submanifolds in part $(i)$.
\vspace{1ex}

In the following global result the assumption that $L^m$ is Einstein 
is replaced by weaker conditions on its scalar curvature.

\begin{theorem}\label{scal}
Let $M^n=L^m\times_\varphi F^{n-m}$, $n\geq m+3\geq 8$, 
be a complete Einstein manifold with Ricci curvature $\rho$. 
Assume that no open subset of $L^m$ admits an isometric 
immersion into $\R^{m+1}$ and that its scalar curvature 
$S_L$ satisfies one of the following:
\begin{itemize}
\item[(i)] $S_L$ is constant,
\item[(ii)] $S_L\leq (2m-n)\rho$.
\end{itemize}
Then any isometric immersion  $f\colon M^n\to\R^{n+2}$ is 
globally a cylinder as in case $(i)$ of Theorem 1.
\end{theorem}

The above results reduce the problem to the classification of the 
Ricci flat submanifolds in Euclidean space with codimension two. 
It is known that these submanifolds always have flat normal bundle; 
see \cite{Co} or Exercise $3.18$ in \cite{DT}. But regardless of 
that simplification, the task seems to be rather difficult due to 
the reason alrady pointed out, namely, the weakness of the Einstein 
assumption in our context. To sense the level of difficulty see 
Corollary $2$ in \cite{VZ}. On the other hand, an abundance of examples 
of Ricci flat submanifolds comes out from Examples $1$ in \cite{DOV2}.

\section{Some general facts}

Throughout the paper $F^m(\varepsilon)$, $m\geq 2$, denotes an 
Einstein manifold with Ricci curvature $(m-1)\varepsilon$, 
$\varepsilon=1,-1,0$. Of course, this is the case of the unit 
sphere $\Sf^m(1)$ and the unit hyperbolic space $\Hy^m(-1)$.

\begin{proposition} 
The warped product $M^n=L^m\times_\varphi F^{n-m}(\varepsilon)$, 
$n-m\geq 2$, is an Einstein manifold of Ricci curvature $\rho$ 
if and only if the warping function $\varphi\in C^\infty(L)$ 
satisfies 
\be\label{eq1}
(n-m)\hess\varphi(X,Y)=(\Ric_L(X,Y)-\rho\<X,Y\>_L)\varphi
\ee
for any vector fileds $X,Y\in\mathfrak{X}(L)$ and that
\be\label{eq2}
\Delta\varphi
+\frac{n-m-1}{\varphi}(\|\n\varphi\|^2
-\varepsilon)+\rho\varphi=0
\ee
where $\hess\varphi$ denotes the Hessian and $\Delta\varphi$ the
Laplacian of $\varphi$. Then $M^n$ and $L^m$ are both Einstein 
manifolds with Ricci curvatures $\rho$ and $\mu$, respectively, 
if and only if $\varphi\in C^\infty(L)$ satisfies  
\be\label{eineq1}
(n-m)\hess\varphi(X,Y)=(\mu-\rho)\varphi \<X,Y\>_L
\ee
for any $X,Y\in\mathfrak{X}(L)$ and
\be\label{eineq2}
\|\n\varphi\|^2
=\varepsilon-\frac{m\mu+(n-2m)\rho}{(n-m)(n-m-1)}\varphi^2\cdot
\ee
\end{proposition}

\proof Equations \eqref{eq1} and \eqref{eq2} are the equations 
9.107b) and 9.107c) in \cite{Be}.\qed

\begin{proposition}\label{coordinates}
Let $M^n=L^m\times_\varphi F^{n-m}(\varepsilon)$ 
and $L^m$,  $n\geq m+2\geq 5$, be Einstein manifolds with Ricci 
curvature $\rho$ and $\mu$, respectively. Assume that 
$\n\varphi\neq 0$ at any point. Then 
$L^m$ is locally a warped product 
$I\times_{\varphi^\prime} N^{m-1}$ where $N^{m-1}$ is an 
Einstein manifold with Ricci curvature 
$\varepsilon\rho(m-2)/(n-1)$, $I\subset\R$ is an open interval 
and $\varphi\in C^\infty(I)$. Moreover, we have:
\begin{itemize}
\item[(i)] If $\rho=0$, then $\varepsilon=1$ and $\varphi(t)=t>0$.
\item[(ii)]  If $\rho>0$, then $\varepsilon=1$ and 
$$
\varphi(t)=a\cos(t\sqrt{\rho/(n-1)})+b\sin(t\sqrt{\rho/(n-1)}
$$ 
where $t\in (0,\pi/2\sqrt{(n-1)/\rho})$ and $a,b>0$ satisfy 
$a^2+b^2=(n-1)/\rho$.
\item[(iii)]  If $\rho<0$, then $\varepsilon=-1,0,1$ and
$$
\varphi(t)=a\cosh(t\sqrt{-\rho/(n-1)})+b\sinh(t\sqrt{-\rho/(n-1)})
$$ 
where $t\in\R$ and $a,b\in\R$ satisfy $a^2-b^2=\varepsilon(n-1)/\rho$.
\end{itemize}
\end{proposition}

\proof  The first part of the statement follows from \eqref{eineq1} 
and Brinkmann's theorem (\cite[Th.\ 4.3.3]{Pet}). Since the
warping function depends only on $t\in I$ 
then $\hess\varphi=\varphi'' \< \,,\,\>_L$ 
and, consequently, \eqref{eineq1} and 
\eqref{eineq2} become 
\be\label{f1}
(n-m)\varphi''=(\mu-\rho)\varphi
\ee
and 
\be\label{f2}
(\varphi')^2
=\varepsilon-\frac{m\mu+(n-2m)\rho}{(n-m)(n-m-1)}\varphi^2\cdot
\ee
Differentiating \eqref{f2} and using \eqref{f1} we obtain
\be\label{rel}
(n-1)\mu=(m-1)\rho.
\ee
Hence \eqref{f1} and \eqref{f2} can be written as
\be\label{f1r}
(n-1)\varphi''+\rho\varphi=0
\ee
and 
\be\label{f2r}
(\varphi')^2=\varepsilon-\frac{\rho}{n-1}\varphi^2.
\ee
If follows from \eqref{rel} and Proposition $9.106$ of 
\cite{Be} that
$$
\Ric_N=\frac{1}{n-1}\left((m-1)\rho(\varphi')^2
+(n-1)\varphi'\varphi'''+(n-1)(m-2)(\varphi'')^2\right)\<\,,\,\>_N.
$$
Hence, taking into account \eqref{f1r} and  \eqref{f2r} we obtain
$$
(n-1)\Ric_N=(m-2)\varepsilon\rho\<\,,\,\>_N.
$$
Finally, solving \eqref{f1r} and using \eqref{f2r} gives the 
expressions for $\varphi$ stated in parts $(i)$ to $(iii)$.
\qed

\section{The proofs}

The proofs of the theorems given in the introduction rely heavily 
on results obtained in \cite{DT0} for isometric immersions of warped 
products into space forms. From there we extract the facts given 
in this section. \vspace{1ex}

In the case of hypersurfaces we have the following:

\begin{proposition}\label{lemma2}
Let $f\colon M^n=L^m\times_\varphi F^{n-m}\to\R^{n+1}$, $m\geq 1$ and 
$n\geq 2$, be an isometric immersion with $M^n$ free of flat points. 
Then $f= \psi\circ(F\times G)$, where
$\psi\colon V^{m+k_1}\times_\sigma\Sf^{n-m+k_2}(r)\to\R^{n+1}$ for
$V^{m+k_1}\subset\R^{m+k_1}$, is a warped product representation of 
$\R^{n+1}$, $k_1+k_2=1$, and the maps $F\colon L^m\to V^{m+k_1}$ 
and $G\colon F^{n-m}\to\Sf^{n-m+k_2}(r)$ are isometric immersions.
\end{proposition}

Recall that an \emph{$(n-m)$-rotational submanifold} 
$f\colon M^n\to\R^{n+2}$, 
$n>m$, with axis $\R^{m+1}$ over a submanifold 
$g\colon L^m\to\R^{m+2}$ is the $n$-dimensional submanifold 
generated by the orbits of the points of $g(L)$ (disjoint 
from $\R^{m+1}$) under the action of the subgroup $SO(n-m+1)$ 
of $SO(n+2)$ which keeps pointwise $\R^{m+1}$ invariant.
Such a submanifold can be parametrized as follows: The manifold 
$M^n$ is isometric to (an open subset of) a warped product 
$L^m\times_\varphi\Sf^{n-m}(1)$ and there is an orthogonal 
splitting $\R^{m+2}=\R^{m+1}\oplus\spa\{e\}$, $\|e\|=1$, such 
that the \emph{profile} $g\colon L^m\to\R^{m+2}$ of $f$ is 
an isometric immersion given by $g=(h,\varphi)$ where 
$h\colon L^m\to\R^{m+1}$ and $\varphi=\<g,e\>>0$. Then,  
we have that $f\colon L^m\times_\varphi\Sf^{n-m}(1)\to\R^{n+2}$
is given by
\be\label{paramet}
f(z,y)=(h(z),\varphi(z)\phi(y))
\ee
where $\phi\colon\Sf^{n-m}(1)\to\R^{n-m+1}$ denotes the 
inclusion $\Sf^{n-m}(1)\subset\R^{n-m+1}$. 
\vspace{1ex}

The following result deals with the case of codimension two.

\begin{proposition}\label{lemma} Let 
$f\colon M^n=L^m\times_\varphi F^{n-m}\to\R^{n+2}, 
n\geq m+3\geq 5$, be an isometric immersion. Assume that no open 
subset of $L^m$ admits an isometric immersion into $\R^{m+1}$.
Then there is an open dense subset of $M^n$ such that each point 
lies in an open product neighborhood 
$U=L^m_0\times F^{n-m}_0\subset M^n$ with $L^m_0\subset L^m$ and 
$F^{n-m}_0\subset F^{n-m}$ such that one of the following holds.
\begin{itemize}
\item[(i)] $f|_U$ is an $(n-m)$-cylinder, that is, there is 
an isometric immersion $h\colon L^m_0\to\R^{m+2}$ and an isometry 
$j$ of $F^{n-m}_0$ with open subset of $\R^{n-m}$ such that 
$f|_U=h\times j$.
\item[(ii)] $f|_U$ is an $(n-m)$-rotational submanifold parametrized 
as in \eqref{paramet}.
\end{itemize}
\end{proposition}

\subsection{The local result}

\noindent {\it Proof of Theorem \ref{local}}:
Up to an homothety, we may assume that $F^{n-m}$ is just 
$F^{n-m}(\varepsilon)$. According to Proposition \ref{lemma}, 
there is an open dense subset of $M^n$ such that each point 
lies in an open product neighborhood $U=L^m_0\times F^{n-m}_0$ 
with $L^m_0\subset L^m$ and $F^{n-m}_0\subset F^{n-m}(\varepsilon)$ 
so that $f|_U$ either is an $(n-m)$-cylinder being $\varphi$ 
constant on $U$ and which consequently is Ricci flat, or 
has to be an $(n-m)$-rotational submanifold. 

Hereafter, we assume that $f|_U$ is an $(n-m)$-rotational submanifold.
First observe that $\varphi$ is not constant on any open subset of 
$L_0^m$ since, otherwise, such a subset would admit an isometric
immersion into $\R^{m+1}$ contradicting the assumption that 
$L^m$ is nontrivial. Thus we may assume that $\n\varphi\neq0$ on $U$. 
Let $\rho$ be the Ricci curvature of $M^n$. From Proposition 
\ref{coordinates} we have that $L^m$ is locally a warped product 
$I\times_{\varphi^\prime} N^{m-1}$ where $N^{m-1}$ is an Einstein 
manifold with Ricci curvature $(m-2)\rho/(n-1)$ and 
$\varphi\in C^\infty(I)$ is given by parts $(i)$ to $(iii)$.

We claim that $M^n$ is Ricci flat. At first we prove that its Ricci
curvature cannot be negative. To the contrary, let us suppose
that $\rho<0$. Then Proposition \ref{coordinates} yields 
$$
\varphi(t)=a\cosh(t\sqrt{-\rho/(n-1)})+b\sinh(t\sqrt{-\rho/(n-1)})
$$
where $t\in I$ and $a,b\in\R$ satisfy $a^2-b^2=(n-1)/\rho$. Setting
$$
a=\sqrt{-(n-1)/\rho}\sinh(t_0\sqrt{-\rho/(n-1)})
$$
and
$$
b=\sqrt{-(n-1)/\rho}\cosh(t_0\sqrt{-\rho/(n-1)})
$$
for some $t_0$, we have
$$
\varphi(t)=\sqrt{-(n-1)/\rho}\sinh((t+t_0)\sqrt{-\rho/(n-1)})
$$
for all $t\in I$. Since we have that the profile
$g\colon L_0^m=I\times_{\varphi'}N^{m-1}\to\R^{m+2}$ of the
$(n-m)$-rotational submanifold  $f|_U$ is
given by $g(t,x)=(h(t,x),\varphi(t))$ where $h\colon L_0^m
=I\times_{\varphi'}N^{m-1}\to \R^{m+1}$, then
$$
1=\<g_*\partial_t,g_*\partial_t\>
=\<h_*\partial_t,h_*\partial_t\>+(\varphi'(t))^2,
$$
which contradicts that $\varphi'(t)>1$ for $t\neq -t_0$.

Next we show that the Ricci curvature cannot be positive. 
To the contrary, if $\rho>0$ then Proposition \ref{coordinates} 
gives that
$$
\varphi(t)=a\cos(t\sqrt{\rho/(n-1)})+b\sin(t\sqrt{\rho/(n-1)})
$$
where $a,b\in\R$ satisfy $a^2+b^2=(n-1)/\rho$.
Setting
$$
a=\sqrt{(n-1)/\rho}\cos(t_0\sqrt{\rho/(n-1)})\;\;\mbox{and}\;\;
b=\sqrt{(n-1)/\rho}\sin(t_0\sqrt{\rho/(n-1)})
$$
for some $t_0$, we obtain 
$$
\varphi(t)=\sqrt{(n-1)/\rho}\cos((t-t_0)\sqrt{\rho/(n-1)})
$$
for $t\in I$. As above the profile of $f|_U$ is
$g(t,x)=(h(t,x),\varphi(t))$ where $h\colon L_0^m
=I\times_{\varphi'}N^{m-1}\to\R^{m+1}$. The metric induced by $h$ is
\begin{align*}
\<\,,\,\>_h
&=(1-(\varphi'(t))^2)dt^2+(\varphi'(t))^2\<\,,\,\>_{N^{m-1}},\\
&=\cos^2\left((t-t_0)\sqrt{\rho/(n-1)})dt^2
+\sin^2((t-t_0)\sqrt{\rho/(n-1)}\right)\<\,,\,\>_{N^{m-1}}.
\end{align*}
Then setting $s=\sqrt{(n-1)/\rho}\sin((t-t_0)\sqrt{\rho/(n-1)})$,
we obtain
$$
\<\,,\,\>_h=ds^2+\frac{\rho s^2}{n-1}\<\,,\,\>_{N^{m-1}}.
$$
Thus we have an isometric immersion
$h\colon J\times_s K^{m-1}\to\R^{m+1}$, where $J\subset \R$
is an open interval and 
$K^{m-1}=(N^{m-1}$, $\rho/(n-1)\<\,,\,\>_{N^{m-1}})$ is 
an Einstein manifold with Ricci curvature $m-2$. Then 
Proposition \ref{lemma2} yields that $h$ is of the 
form $h=\psi\circ(f_0\times g_0)$, where 
$\psi\colon V^{k_1+1}\times_\sigma\Sf^{m+k_2-1}(r)\to\R^{m+1}$ 
is a warped product representation of $\R^{m+1}$  such that 
$k_1+k_2=1$,
the maps $f_0\colon J\to V^{k_1+1}\subset\R^{k_1+1}$ and
$g_0\colon K^{m-1}\to\Sf^{m+k_2-1}(r)$ are isometric
immersions and $s=\sigma\circ f_0$. Recall that  $\psi$ 
is the map defined by
$$
\psi(p_0,p_1)=p_0+\sigma(p_0)(p_1-q),
$$
where $q$ is a point in the intersection $V^{k_1+1}\cap \Sf^{m+k_2-1}(r)$ 
and $\sigma\colon V^{k_1+1}\to\R_+$ is the linear function given by
$\sigma(p_0)=\<p_0,q\>/r^2$.

According to the above we need to distinguish two cases. At first 
let $k_1=1$. Hence $f_0$ is just a curve in $V^2\subset\R^2$ and 
$g_0$ an isometry of $K^{m-1}$ to an open subset of $\Sf^{m-1}(r)$. 
Moreover, since $K^{m-1}$ is an Einstein manifold with Ricci 
curvature $m-2$ then $r=1$. Hence the  metric induced by $h$ is 
the flat metric
$$
\<\,,\,\>_h=ds^2+s^2\<\,,\,\>_{\Sf^{m-1}(1)}.
$$
For simplicity we choose $q=(1,0,\ldots,0)$ and let 
$f_0(s)=(x(s),0,\ldots,0,y(s))$ be a parametrization.
Clearly $f_0(s)$ is a unit speed curve. Since $s=\sigma(f_0(s))$
then $x(s)=s$ and hence $y(s)$ is constant. Thus $h(L_0^m)$ lies
in a hyperplane of $\R^{m+1}$ and therefore the profile in a 
hyperplane of $\R^{m+2}$, contradicting the assumption that
the Einstein manifold $L^m$ is nontrivial.

Now let $k_2=1$. Then $f_0$ is the identity and $g_0$ an Einstein 
hypersurface in $\Sf^m(r)$. For simplicity choosing
$q=(r,0,\dots,0)$ in the warped product representation 
of $\R^{m+1}$, then $h\colon J\times_sK^{m-1}\to\R^{m+1}$ 
is the cone given by
$$
h(s,z)=\frac{s}{r}g_0(z).
$$
This implies that the induced metric by $h$ is
$$
\<\,,\,\>_h=ds^2+(s/r)^2\<\,,\,\>_{K^{m-1}}.
$$
On the other hand, since we had that
$$
\<\,,\,\>_h=ds^2+s^2\<\,,\,\>_{K^{m-1}},
$$
hence $r=1$ and $h\colon I\times_{\varphi'}N^{m-1}\to\R^{m+1}$
is given by
$$
h(t,z)=\sqrt{(n-1)/\rho}\sin((t-t_0)\sqrt{\rho/(n-1)})g_0(z).
$$
It follows from the parametrization \eqref{paramet} that
the rotational submanifold is contained in the sphere
$\Sf^{n+1}(\sqrt{(n-1)/\rho})$ and hence has constant sectional
curvature according to the classification of Einstein spherical 
submanifolds due to Fialkow \cite{Fi}. But this contradicts our 
assumption on $M^n$ and proves the claim.
\vspace{1ex}

Since $M^n$ is Ricci flat then  Proposition 
\ref{coordinates} yields $\varphi(t,x)=t$. 
The profile $g\colon L_0^m=I\times N^{m-1}\to\R^{m+2}$ 
of the $(n-m)$-rotational submanifold  $f|_U$ is
$g(t,x)=(h(t,x),t)$ where 
$h\colon L_0^m=I\times N^{m-1}\to \R^{m+1}$. That
$$
1=\<g_*\partial_t,g_*\partial_t\>
=\<h_*\partial_t,h_*\partial_t\>+1,
$$
implies that $h$ is independent of $t$. Hence $g=f_0\times id_I$ 
where $f_0=h$ is a nonflat Ricci flat submanifold. Now it follows 
from \eqref{paramet} that the submanifold is as in $(ii)$. 
To conclude the proof, observe that $\|\n\varphi\|$ has to take 
the same constant value in all of $M^n$.\qed

\begin{remarks} {\em
\noindent $(1)$ The first assumption  in Theorem \ref{local} that $m>4$
is necessary. In fact, if in part $(ii)$ we have $m=4$ then $N^3$ 
being Einstein has constant sectional curvature and the same would 
be the case of $L^4_0$, but that possibility has been excluded.
\vspace{1ex}

\noindent $(2)$  The first assumption on $L^m$ in Theorem \ref{local} 
cannot be dropped as shown by the following example of an Einstein
manifolds that admits an isometric immersion in codimension two 
but not as a cylinder.
\vspace{1ex}

\noindent 
Let $M^n=\Sf^m(\sqrt{(m-1)/\rho})\times\Sf^{n-m}(\sqrt{(n-m-1)/\rho})$ 
and let $f\colon M^n\to\R^{n+2}$ be the product of the inclusions 
$\Sf^{n_j}(r_j)\subset\R^{n_j+1}$, where $n_j=m,n-m$ and  
$r_j=\sqrt{(n_j-1)/\rho}$. Notice that this example justifies the 
necessity of the assumption on $L^m$ in Theorem \ref{scal}.
}\end{remarks}

\subsection{The global result}

\noindent {\it Proof of Theorem \ref{scal}}:  
We assume that $F^{n-m}=F^{n-m}(\varepsilon)$ and 
by Proposition \ref{lemma} there is an open dense 
subset of $M^n$ so that each point lies in an open product 
neighborhood $U=L^m_0\times F^{n-m}_0\subset M^n$ with 
$L^m_0\subset L^m$ and $F^{n-m}_0\subset F^{n-m}(\varepsilon)$ 
such that $f|_{U}$ is either an $(n-m)$-cylinder and thus 
$\varphi$ is constant on $L_0^m$, or it is an $(n-m)$-rotational 
submanifold. In the latter case, $f|_{U}$ is parametrized by 
\eqref{paramet} and hence $\varphi=\<g,e\>$, where 
$g\colon L^m_0\to \R^{m+2}$ is the profile and $e$ is a unit 
vector in $\R^{m+2}$. Since $\n\varphi$ is the tangent 
component of $e$ then $\|\nabla\varphi\|\leq 1$. Therefere, 
we have in either case that $\|\nabla\varphi\|\leq 1$ on $L_0^m$, 
and by continuity on all of $L^m$. 

Taking traces in \eqref{eq1} gives
\be\label{eq32}
(n-m)\Delta \varphi=\varphi (S_L-m\rho).
\ee
Then combining with \eqref{eq2} yields
\be\label{eq42}
S_L+(n-2m)\rho
=(n-m)(n-m-1)\frac{\varepsilon-\|\n \varphi\|^2}{\varphi^2}\cdot
\ee

We claim that $\varphi$ is constant on $L^m$. Suppose to the
contrary that $\varphi$ is not constant. It follows from
Proposition \ref{lemma} that $M^n$ has an open subset $V$ such 
that $f|_V$ is an $(n-m)$-rotational submanifold. Thus 
$\varepsilon=1$ and \eqref{eq42} becomes
\be\label{eq52}
\|\n \varphi\|^2=1-\frac{S_L+(n-2m)\rho}{(n-m)(n-m-1)}\varphi^2\cdot
\ee
Since $\|\nabla\varphi\|\leq 1$, we have $S_L\geq(2m-n)\rho$. Clearly, 
we have to distinguish two cases.

Suppose that $S_L=(2m-n)\rho$. Then \eqref{eq52} implies that
$\|\nabla\varphi\|= 1$ everywhere and hence the integral curves of 
$\n\varphi$ are geodesics. If $\gamma(t)$, $t\in \R,$ is such
a curve then $(\varphi\circ\gamma)'(t)=1$. Hence 
$(\varphi\circ\gamma)(t)=t+c_0$ for any $t\in\R$, which is a 
contradiction.

Suppose now that $S_L$ is constant and $S_L>(2m-n)\rho$. 
It follows from  \eqref{eq52} that $\varphi\leq a$, where
$$
a=\sqrt{\frac{(n-m)(n-m-1)}{S_L+(n-2m)\rho}}.
$$
In addition,  \eqref{eq52} implies that the set $M_c$ of the 
critical points of the warping function coincides with the 
set of points where it attains the value $a$. We consider on 
$M\smallsetminus M_c$ the unit vector field
$v=\nabla\varphi/\|\nabla\varphi\|$. Using \eqref{eq52} and 
$$
2\n_{\n\varphi}\n\varphi=\n\|\n\varphi\|^2,
$$
where $\n$ stands for the Levi-Civita connection on $L^m$, we 
obtain that the integral curves of $v$ are unit speed geodesics. 
Let $\gamma(s)$ be such a curve and consider the function 
$u=\varphi\circ \gamma$. Then \eqref{eq52} yields
\be\label{eq62}
(u'(s))^2=1-\frac{1}{a^2}(u(s))^2
\ee
for all $s\in\R\smallsetminus A$, where 
$A=\{s\in\R:\gamma(s)\in M_c\}$.  But besides the constant solution
we have that $u(s)=a\sin(c\pm s/a)$, $c\in\R$ which is defined on 
the whole real line and then has to be  discarded since it takes 
negative values. Thus the claim has been proved.

Since $\varphi$ is constant then $\varepsilon=0$. In fact, if 
$\varepsilon=1$ we have by Proposition \ref{lemma} that $f$ 
is locally $(n-m)$-rotational parametrized by \eqref{paramet}, thus
contradicting the assumption that no open subset of $L^m$ admits 
an isometric immersion into $\R^{m+1}$. Hence, from \eqref{eq42} 
and that $\varphi$ is constant, we obtain $S_L=(2m-n)\rho$.
But since \eqref{eq32} yields that $S_L=m\rho$,  then $M^n$ is 
Ricci flat. 

We have that there in no open subset $V\subset M^n$ 
such that $f|_V$ is an $(n-m)$-rotational submanifold 
parametrized by \eqref{paramet}. In fact, if otherwise the 
profile of $f|_V$ would lie in a hyperplane of $\R^{m+2}$, 
and consequently, an open subset of $L^m$ would admit an 
isometric immersion into $\R^{m+1}$. Hence $f$ is locally 
an $(n-m)$-cylinder, and since $\rho=0$ then $L^m$ is 
Ricci flat. That $f$ is globally a cylinder over a complete 
Ricci flat submanifold follows from Theorem \ref{local} and 
the completeness of $M^n$.\qed

\section*{Acknowledgment}

The first and third authors thank the Mathematics Department of 
the University of Murcia where part of this work was developed for 
the kind hospitality during their visits.
\medskip

This research is part of the grant PID2021-124157NB-I00, funded by\\
MCIN/ AEI/10.13039/501100011033/ ``ERDF A way of making Europe"

\noindent Marcos Dajczer\\
IMPA -- Estrada Dona Castorina, 110\\
22460--320, Rio de Janeiro -- Brazil\\
e-mail: marcos@impa.br
\bigskip

\noindent Christos-Raent Onti\\
Department of Mathematics and Statistics\\
University of Cyprus\\
1678, Nicosia -- Cyprus\\
e-mail: onti.christos-raent@ucy.ac.cy
\bigskip

\noindent Theodoros Vlachos\\
University of Ioannina \\
Department of Mathematics\\
Ioannina -- Greece\\
e-mail: tvlachos@uoi.gr

\end{document}